# Congestion Management by a Real-Time Optimal Dispatch through Balancing Mechanism Method for Electricity Markets

Barati Masoud, Shayanfar Heidar Ali, Kazemi Ahad

*Abstract*— A real-time optimal dispatch method for unbundled electricity markets is proposed in this paper. With this method, pool energy auction market, ancillary services market, and bilateral contract market can be coordinated by the ISO through a balancing market for the purpose of system security. To meet real-time imbalance and to eliminate transmission congestion, operating reserves, curtailment of bilateral contracts, and supplemental energy bids in the real-time balancing market (REM), can all be called upon in light of their bidding prices. The replacement of the operating reserves used to provide balancing energy is incorporated into the objective of the real-time optimal dispatch problem. A modified P-Q decoupled optimal power flow (OPF) is employed to implement this method. Spot pricing for the real-time dispatch is presented. The IEEE 30-bus test system is studied to illustrate the proposed framework.

**Index Terms**-Ancillary services, balancing market, bilateral contract, congestion management, coordinated dispatch, optimal power flow (OPF), deregulation, energy market.

## I. INTRODUCTION

However, the problem still remaining for an ISO to resolve is how to use all the possible resources during the real-time execution of various electricity commodity contracts efficiently and coordinately to ensure the system security. The main difficulties occurring in this real-time coordinated dispatch problem could be:
- how to dispatch the agreed system reserves contracts together with the supplemental energy bids in the balancing market;
- with the trend that more and more bilateral contracts are used to trade electricity, how to eliminate network congestion if the resources in the balancing market are not enough;
- To maintain the system security level, how to purchase replacement operating reserves if any of the pre-arranged operating reserves are called upon to provide energy for real-time system balancing or congestion management [2].

To resolve these difficulties, a new framework for real-time dispatch of unbundled electricity markets is proposed in this article. Under this framework, almost all the contracts in various markets can be dispatched and coordinated by submitting their adjustment bids to the balancing market. In particular, some bilateral contracts can be adjusted if the congestion of the network is very serious [3, 4]. Demand side participants are encouraged to play an active role in the competition of the real-time balancing market. A modified P-Q decoupled OPF is applied to solve this problem. The objective of the P sub-problem is to coordinate the dispatch among the pool auction contracts, the bilateral contracts and some operating reserve contracts, in accordance with the cost which is determined by the bids submitted by these contracts to the balancing market. The spot pricing and the two possible settlement methodologies are analyzed as well. A 5-bus test system and the IEEE 30-bus test system are studied to illustrate the proposed framework and its mathematical solution.

## II. FRAMEWORK OF REAL-TIME COORDINATED DISPATCH

Real-time Balancing Market and Coordinated Dispatch (RBM) play a key role in real-time operation of electricity markets. In RBM the ISO is responsible for dispatching all the available resources to meet imbalances between actual and scheduled load and generation and alleviate network congestion. The ISO will select the least cost resources to meet these imbalances. In addition, the ISO may also need to purchase replacement ancillary services if any services arranged in advance are used to provide balancing energy [2]. The purpose of RBM is to establish a fully open market-based mechanism for all the market participants to take part in the real-time competition. In RBM, all the generators and consumers can submit to the ISO their incremental and decremental bids for providing balancing energy and their capacity bids for the replacement of operating reserves which are used as balancing energy. In the NETA of the UK market, these incremental/decremental bids are called "pairs of offer and bid". Because the ISO does not know the price information of bilateral contracts and the modification of a bilateral contract could involve both sides of a contract, it is very difficult to find a proper way to change them for the purpose of reducing transmission congestion. One method is that both parties submit their own supplemental bids in RBM separately, just like the other participants in the pool. The ISO does not take into account the content of the bilateral contract during settlement. However, in the event of a high percentage of bilateral contracts not producing enough voluntary supplemental bids from them or both parties of a bilateral contract are willing to be curtailed by the same amount to simplify the settlement, other methods are needed. In the proposed framework, every bilateral contract should submit a sort of compensative price that both parties of the contract are willing to accept if the curtailment needs to be

The authors are with the Center of Excellent for Power System Automation and Operation Department of Electrical Engineering, Iran University of Science & Technology, Narmak 16846-13114 Tehran, IRAN (barati196@gmail.com; hshayanfar@yahoo.com; kazemi@iust.ac.ir).

imposed by the ISO during congestion periods. In accordance with such information, the ISO can reduce the scheduled bilateral contract if the other available resources in RBM are not enough to eliminate the congestion.

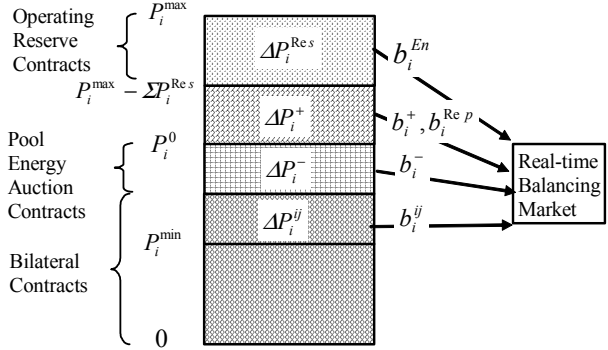

Fig.1. Typical contracts associated with a generator at bus $i$ and their adjustment bids to RBM

In Fig 1,

| | |
|---|---|
| $i,j$ | Index of network buses. |
| $P_i^0$ | Scheduled MW generation of generator or MW load of consumer at bus $i$. |
| $P_i^{ij}$ | Total MW amount of bilateral contracts between buses $i$ and $j$. |
| $P_i^{min}, P_i^{max}$ | Low and high limit of MW at bus $i$ either for generator or for consumer. |
| $P_i^{Res}$ | Operating reserves procured by ISO from participant at bus $i$ in ancillary services market. |
| $\Delta P_i^+, \Delta P_i^-$ | Incremental and decremental MW changes of participant at bus $i$ in real-time balancing market. |
| $\Delta P_i^{ij}$ | Curtailment of bilateral contract between buses $i$ and $j$. |
| $\Delta P_i^{Res}$ | Operating reserves called upon by ISO in the real-time balancing market out of $P_i^{Res}$. |
| $\Delta P_i^{Rep}$ | Reserves at bus $i$ procured by ISO in real-time balancing market to replace the used operating reserves. |
| $b_i^+, b_i^-$ | Incremental and decremental bidding prices of participant at bus $i$ in RBM ($/MW). |
| $b_i^{ij}$ | Compensative price for the curtailment of bilateral contract between buses i and j ($/MW). |
| $b_i^{En}$ | Bidding price of the energy output from $P_i^{Res}$ ($/MW). |
| $b_i^{Rep}$ | Bidding price to provide $\Delta P_i^{Rep}$ at bus $i$ ($/MW). |

Some typical contracts of active power associated with a generator at bus $i$ and their bids submitted to the ISO in RBM are shown in Fig.1, in which the adjustment of generation at bus $i$ is divided into four separate parts based on three unbundled electricity contracts.

### III. MATHEMATICAL MODEL OF THE PROPOSED FRAMEWORK

A modified P-Q decoupled OPF is applied to implement the proposed framework. Both P- and Q sub-problems are presented, and the P sub-problem is analyzed in more detail. Before the trade period, the schedules of unbundled markets should have been established. The principle of the balancing mechanism is to follow the schedules as close as possible and to minimize the cost of real-time dispatch. After linearization at the scheduled operating point, the model of real-time coordinated optimal dispatch can be written in following sections.

### IV. P SUB-PROBLEM

#### A. Objective

To make the settlement of unbundled markets clear and simple, the objective of real-time active power dispatch can be decomposed into four parts.

$$Min \quad C_p = C_{p1} + C_{p2} + C_{p3} + C_{p4} \quad (1)$$

Where

$$C_{p1} = \sum_{i \in PEAM}(b_i^+ \Delta P_i^+ + b_i^- \Delta P_i^-) \quad (2)$$

The total cost for the adjustment of PEAM contracts whose supplemental bids are accepted in RBM. But for accepted decremental bids, the participants can not get any payment for this part of reduced energy from the ISO in PEAM settlement;

$$C_{p2} = \sum_{i \in PAAM} b_i^{En} \Delta P_i^{Res} \quad (3)$$

The cost for calling upon energy from the agreed operating reserve contracts in RBM;

$$C_{p3} = \sum_{i \in BCM} \sum_{\substack{j \in BCM \\ j<i}}(b_i^{ij} \Delta P_i^{ij}) \quad (4)$$

is the cost for curtailment of bilateral contracts in RBM. This curtailment should be done at the both sides of bilateral contract; and

$$C_{p4} = \sum_{i \in RBM}(b_i^{Rep} \Delta P_i^{Rep}) \quad (5)$$

is the cost of reserves re-procured by ISO in RBM to replace the used operating reserves. In the event of an unscheduled increase in system load, the ISO will be required to purchase additional reserves. In addition, if any of the pre-arranged operating reserves are used to meet a real-time imbalance, the ISO will be required to purchase replacement operating reserves. The reserve market and the balancing energy market will be still settled separately. Embedding the cost of replacement of operating reserves in the objective function of real-time dispatch is to find a global optimal solution. It can be regarded as the real-time joint dispatch of energy and reserves.

Obviously, the objective is to minimize the modification on all the scheduled contracts in light of the associated dispatch cost. $\Delta P_i^+$, $\Delta P_i^-$, $\Delta P_i^{ij}$, $\Delta P_i^{Res}$ and $\Delta P_i^{Rep}$ are treated as independent control variables during the optimization process. But their upper/lower limits are coupled with each other. All the bidding curves could be multi-step. The incremental bidding price is higher than decremental bidding price while the curtailment price of bilateral contracts is much higher than the other two. The reason is that increasing output needs more fuel cost and the curtailment to a bilateral contract will affect the financial interests of both parties.

*B. Equality Constraints*

The equation below is the nodal active power flow balance equation of bus $i$:

$$(-1)^\beta [\Delta P_i^+ - \Delta P_i^- - \sum_{\substack{j \in BCM \\ j \neq i}} (\Delta P_i^{ij})] + \Delta P_i^{Res} \quad (6)$$
$$+ \sum_{\substack{j=1 \\ j \neq s}}^{n} B'_{ij} \Delta \theta_j - \Delta P_i^{Loss} = 0$$

$\beta = 0, \text{if } i \in G; \beta = 1, \text{if } i \in C$.

Where $G$ is the set of buses connected with generators, $C$ is the set of buses connected with consumers, $n$ is the total number of network buses, $s$ is the index of the slack bus, $B'_{ij}$ is an element of matrix $\mathbf{B}'$ which is the inverse reactance of branch $ij$. $\Delta P_i^{Loss}$ Is the summed change of losses on branches that are connected to bus $i$ and flows on the branches are flowing to bus $i$. Using the piece-wise linear loss model.

The requirement of re-procurement to replace used operating reserves is given by

$$\sum_{i \in RBM} \Delta P_i^{Rep} = \alpha \sum_{i \in PAAM} \Delta P_i^{Res} \quad (7)$$

Where $\alpha$ is decided by the ISO according to how many new operating reserves are procured to compensate for the used reserves in real-time dispatch. $0 \leq \alpha \leq 1$.

*C. Inequality Constraints*

The changing ranges of the various control variables are given by

$$0 \leq \Delta P_i^+ + \Delta P_i^{Rep} \leq \Delta P_i^{+,max} = P_i^{max} - P_i^0 - P_i^{Res} \quad (8)$$

$$0 \leq \Delta P_i^- \leq \Delta P_i^{-,max} = P_i^0 - \max(\sum_{\substack{j=1 \\ j \neq n}}^{n} P_i^{ij}, P_i^{min}) \quad (9)$$

$$0 \leq \Delta P_i^{Res} \leq P_i^{Res} \quad (10)$$

$$0 \leq \Delta P_i^{ij} \leq P_i^{ij}, 0 \leq \sum_{\substack{j=1 \\ j \neq i}}^{n} \Delta P_i^{ij} \leq \min(\sum_{\substack{j=1 \\ j \neq i}}^{n} P_i^{ij}, P_i^{min}) \quad (11)$$

$$-P_l^{max} - P_l^0 = \Delta P_l^{min} \leq \Delta P_l \leq \Delta P_l^{max} = P_l^{max} - P_l^0 \quad (12)$$

The constraint for real power flow change on branch $l$.

*D. Pricing for Real-Time Active Power Dispatch*

From (6), we can have the system active power balance equation as:

$$\sum_{i=1}^{n} \Delta P_i - \Delta P^{Loss} = 0 \quad (13)$$

Where $\Delta P_i$ is the total change of active power at bus $i$ and $\Delta P^{Loss}$ is the total change of active power losses. $\Delta P^{Loss} = \sum_{i=1}^{n} \Delta P_i^{Loss}$.

Using the theory of spot pricing [1, 5], the LMP of a participant (either a generator or a consumer) at bus $i$ is:

$$\rho_i^P = \lambda - \frac{\partial \Delta P^{Loss}}{\partial \Delta P_i} \lambda - \sum_{l \in B}(\mu_l^{max} - \mu_l^{min}) \frac{\partial \Delta P_l}{\partial \Delta P_i} \quad (14)$$

$$\rho_i^P = \begin{cases} (-1)^\beta (\frac{\partial C_p}{\partial \Delta P_i^+} + \mu_i^+); \text{ or} \\ (-1)^{\beta+1}(\frac{\partial C_p}{\partial \Delta P_i^-} + \mu_i^-); \text{ or} \\ \frac{\partial C_p}{\partial \Delta P_i^{Res}} + \alpha \mu^{Rep} + \mu_i^{Res}. \end{cases} \quad (15)$$

Shown in equation (14), the real-time spot price at bus $i$ can be decomposed into three parts: the system lambda, the active power losses and the congestion management cost.

From equation (15), we can have different forms of $\rho_i^P$ at bus $i$ at which the participant has made some contribution to the real-time dispatch.

*E. Meeting Real-Time Imbalance of Market under Normal Operating Conditions*

Providing there is no serious contingency, the operating reserves and the supplemental energy in RBM should be enough to follow the system load fluctuation. The curtailment of bilateral contracts and load shedding are not necessary in this case. To model this problem, the objective in (1) should be rewritten as:

$$Min\ C_p = C_{p1} + C_{p2} + C_{p4} \quad (16)$$

And all the constraints of $\Delta P_i^{ij}$ should be removed.

Given the system load fluctuation $\Delta P^{sys}$, the bus load change can be expressed as:

$$BL_i = \eta_i \Delta P^{sys} \quad (17)$$

Where $\eta_i$ is a bus load allocation factor, $\sum_{i \in C} \eta_i = 1$.

To meet this system imbalance, the right hand side of equation (6) changes to $BL_i$ from 0. Considering that rapid response generation units are enough to meet the normal system imbalance, to reduce the number of control variables, the nodal active power balance equations (6) for pure load buses can be rewritten as:

$$\sum_{\substack{j=1 \\ j \neq s}}^{n} B'_{ij}\Delta\theta_j - \Delta P_i^{Loss} = BL_i \quad (18)$$

And the system balance equation (13) can be changed to:

$$\sum_{i=1}^{n}\Delta P_i - \Delta P^{Loss} = \Delta P^{sys} \quad (19)$$

Without any branch limit violation, the LMP at bus $i$ is:

$$\rho_i^P = \lambda - \frac{\partial \Delta P^{Loss}}{\partial \Delta P_i}\lambda \quad (20)$$

### F. Replacement of Operating Reserves

From equation (19), the marginal price of the replacement of operating reserves is:

$$\rho^{\mathrm{Re}\,p} = \mu^{\mathrm{Re}\,p} + \mu_i^+ \quad (21)$$

The cost of replacement of operating reserves affects the LMPs through the item $\alpha\mu^{\mathrm{Re}\,p}$ instead of appearing in the equation (14) or (20) directly.

Calling up energy from operating reserves during real-time dispatch, the ISO must pay for not only the energy but also the procurement of the reserve replacement. Therefore, if the supplemental energy in RBM is fast and cheap enough, the reserve contracts signed in PAAM could be left unchanged. As a result, the real-time balancing mechanism expands the concept of operating reserves and gives market participants more competitive opportunities.

### G. Curtailment of Bilateral Contracts

According to the characteristics of bilateral contract, the curtailment should be done to both parties of the contract as has been discussed above.

But sometimes the load in a bilateral contract is too important to be curtailed or the response of the generation unit in a bilateral contract is not fast enough to decrease its output to mitigate the network congestion. In these two cases, the energy imbalance caused by single side curtailment will be taken by the other available resources in RBM. Here, the partner of a bilateral contract without curtailment will not get the compensation payment from the ISO. On the contrary, this participant must pay for the dispatch cost for this imbalance. Then the equation (23) has another form as:

$$\rho_i^P = (-1)^{\beta+1}(\frac{\partial C_p}{\partial \Delta P_i^{ij}} + \mu_i^{ij}) \quad (22)$$

## V. Q SUB-PROBLEM

The Q Sub-problem of real-time coordinated optimal dispatch can be formulated as:

$$\text{Min } \sum_{i=1}^{n}(w_i^+\Delta Q_i^+ + w_i^-\Delta Q_i^-) + \sum_{k=1}^{T}r_k(\Delta t_k^+ + \Delta t_k^-) \quad (23)$$

Subject to:

$$\Delta Q_i^+ - \Delta Q_i^- + \sum_{j=1}^{n} B''_{ij}\Delta V_j + \sum_{k=1}^{T}\partial Q_i/\partial t_k (\Delta t_k^+ - \Delta t_k^-) = 0 \quad (24)$$

$$0 \leq \Delta Q_i^+ \leq Q_i^{\max} - Q_i^0, 0 \leq \Delta Q_i^- \leq Q_i^0 - Q_i^{\min} \quad (25)$$

$$0 \leq \Delta t_k^+ \leq t_k^{\max} - t_k^0, 0 \leq \Delta t_k^- \leq t_k^0 - t_k^{\min} \quad (26)$$

$$V_i^{\min} - V_i^0 \leq \Delta V_i \leq V_i^{\max} - V_i^0 \quad (27)$$

where $w_i^+$ and $w_i^-$ are the reactive power incremental and decremental bidding prices of participant $i$ in RBM respectively, $\Delta Q_i^+$ and $\Delta Q_i^-$ are its increasing output and decreasing output of reactive power respectively, $Q_i^0$, $Q_i^{\min}$, $Q_i^{\max}$ are its current reactive power output or load, minimum reactive power and maximum reactive power respectively. $T$ is the number of transformers, whose tap positions are adjustable, in the system. $\Delta t_k^+$ and $\Delta t_k^-$ are the increasing tap ratio and the decreasing tap ratio of transformer $k$. $r_k$ is the bidding price for the adjustment of tap positions of transformer $k$ if there are any independent transmission companies and the transmission sector is also competitive. $t_k^0, t_k^{\min}, t_k^{\max}$ are the current tap ratios, the minimum and maximum tap ratios of transformer $k$, respectively. $\Delta V_i$ is the change of the voltage magnitude at bus $i$, $V_i^0$, $V_i^{\min}$ and $V_i^{\max}$ are its current value, lower limit and upper limit. $B''_{ij}$ is an element of the matrix $\mathbf{B''}$.

## VI. IMPLEMENTATION

To solve the problem fast and reliably, an AC Power Flow and a Primal-dual Interior Point (IP) Linear Programming (LP) are used to implement the above OPF algorithm. The procedure of the decoupled OPF implemented in this article is as follows:

Step1: Run the AC Power Flow to get the initial state of power system;

Step2: Compute the necessary sensitivities and linearize the constraints;

Step3: Select the control variables of LP according to the bids in RBM. Decide if the bilateral contract curtailment is needed and which curtailment strategy should be used;

Step4: Run LP to solve P sub-problem;

Step5: Run LP to solve Q sub-problem;

Step6: Correct the control variables, and then run AC Power Flow to get the new state of power system;

Step7: Check if all the constraints have been satisfied. If yes, continue; if no, go to step2;

Step8: Obtain the optimal coordinated dispatch strategy.

## VII. SIMULATION

Two test systems are studied to illustrate the proposed method. The IEEE 30-bus test system is studied to show how the real-time coordinated optimal dispatch works and the LMPs obtained from it. Finally, a computational comparison between the Simplex LP method and Interior Point LP method is performed to show the efficiency of the latter. The IEEE 30-bus system is used here to illustrate the proposed coordinated dispatch method. The various bids of participants, including generators, consumers and bilateral contracts, are given in Table I,

TABLE I THE BIDS OF PARTICIPANTS IN RBM

| Partici-pants | Supplemental bids | | Operating Reserves | | | Base Point (MW) | MAX MW | MIN MW |
|---|---|---|---|---|---|---|---|---|
| | Incr. Bids ($/MWh) | Decr. Bids ($/MWh) | Cap. Bids ($/MW) | En. Bids ($/MWh) | Amount (MW) | | | |
| G-1 | 35 | 15 | 3.5 | 35 | / | 138.53 | 200 | 50 |
| G-2 | 15 | 8 | 2.5 | 15 | / | 57.56 | 100 | 20 |
| G-5 | 15 | 8 | 1.5 | 15 | / | 24.56 | 100 | 10 |
| G-8 | 30 | 12 | 1.5 | 15 | 30 | 35.00 | 65 | 10 |
| G-11 | 25 | 10 | 2.5 | 25 | / | 17.93 | 50 | 10 |
| G-13 | 15 | 5 | 1.5 | 15 | / | 16.91 | 50 | 5 |
| C-24 | / | 40 | / | / | / | 8.70 | 15 | 3 |

| Bilateral contract | From | To | Contract Amount (MW) | Curtailment Bids ($/MW) |
|---|---|---|---|---|
| B1 | G-13 | C-30 | 10.6 | 50 |

In which G-1 means generator at bus 1 and C-24 means consumer at bus 24.

*A. Coordinated Dispatch without Network Congestion*

Assume that there is a 100MW increase on system load, which has been distributed to individual buses according to their current load shares. Set $\alpha = 1$, then the obtained optimal dispatch strategy to meet this load fluctuation is shown in Table II, Fig. 2 reveals the two components of LMPs under normal operating conditions, which are system lambda and network losses.

TABLE II THE OPTIMAL DISPATCH STRATEGY TO MEET THE LOAD FLUCTUATION

| Partici-pants | PEAM contracts | | Calling Upon PAAM contracts (MW) | Replacement of Operating Reserves (MW) |
|---|---|---|---|---|
| | Increase (MW) | Decrease (MW) | | |
| G-5 | 75.44 | 0 | 0 | 0 |
| G-8 | 0 | 0 | 24.72 | 0 |
| G-13 | 0 | 0 | 0 | 24.72 |

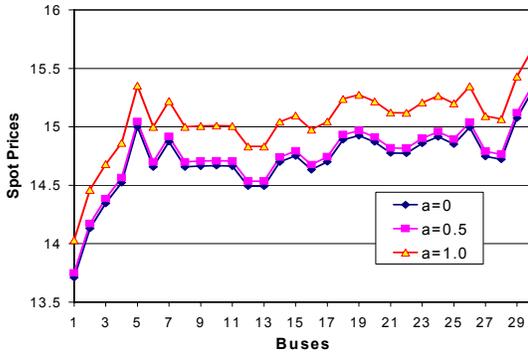

Fig.2. Components of spot prices under normal operating condition

eliminate this light congestion with meeting the given load change is shown in Table III.

In this case, three components of the LMP given in Equation (14), which are system lambda, network losses and congestion management are illustrated in Fig. 4.

The cost of replacement of used operating reserves is not a component of LMPs in Equation (28). However, it can affect the value of system lambda. Fig.7 demonstrates the change of LMPs under different replacement Procurement of used operating reserves by giving $\alpha = 0$, $\alpha = 0.5$, $\alpha = 1$.

*B. Coordinated Dispatch with Network Congestion*

Case 1: Based on the same load fluctuation and the same operating reserves replacement level ($\alpha = 1$) as those in last section, the active power flow on line 36 (from bus-28 to bus-27) is 24.79MW. Now reduce the limit of active power flow on this line to 24MW, the optimal dispatch strategy to

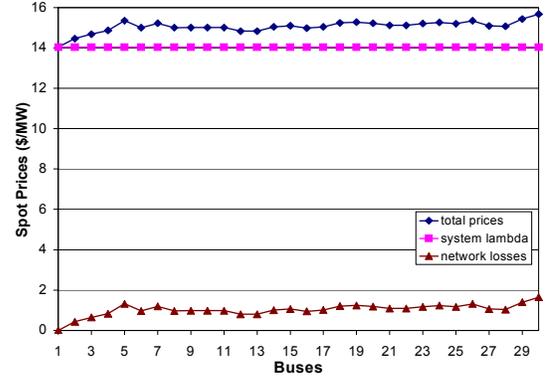

Fig.3. Effects of different replacement procurement of operating reserves on LMPs

TABLE III THE OPTIMAL DISPATCH STRATEGY WITH ACTIVE POWER FLOW VIOLATION ON LINE 36.

| Partici-pants | PEAM contracts | | Calling Upon PAAM contracts (MW) | Replacement of Operating Reserves (MW) |
|---|---|---|---|---|
| | Increase (MW) | Decrease (MW) | | |
| G-5 | 75.44 | 0 | 0 | 0 |
| G-8 | 0 | 0 | 13.22 | 0 |
| G-13 | 11.63 | 0 | 0 | 13.22 |

Compared with Fig. 2, the spot prices in Fig. 4 have not changed a lot at most buses except buses 25-27, 29-30.

That means these several buses have higher sensitivity to the congested line 36 than the other buses. In other word, the consumers at these buses should pay most of the cost caused by network congestion. Because the congestion is very slight, the associated Lagrangian multiplier $\mu_{36}^{\max} = -2.21$ is not terribly big.

**Case 2**: Based on the same operating condition as in case 1, but the congestion branch in this case is line 18. (From bus-12 to bus-15), whose base active power flow is 22.67MW. Now reduce the limit of active power flow on this line to 19MW because of an unexpected contingency, the optimal dispatch strategy to eliminate this congestion with meeting the given load change is shown in Table IV. From Table 4 it can be noticed that this network congestion is so serious that the decremental bids of G-13 is still not enough to eliminate it. Some transaction in bilateral contract B1 has to be curtailed. The cost of this curtailment can not be embedded into the defined spot prices. It must be allocated to market participants separately.

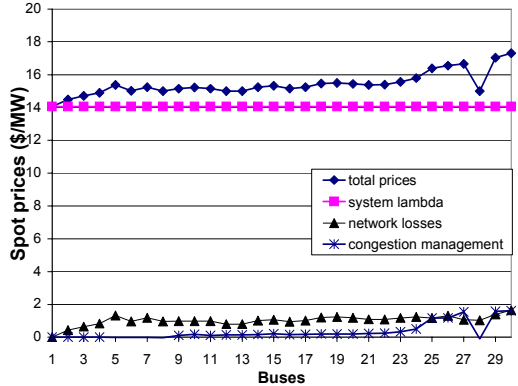

Fig.4. Components of spot prices with network congestion

TABLE IV THE OPTIMAL DISPATCH STRATEGY WITH ACTIVE POWER FLOW VIOLATION ON LINE 18.

| Partici-pants | PEAM contracts | | Calling Upon PAAM contracts (MW) | Curtailment of Bilateral Contracts (MW) | Replacement of Operating Reserves (MW) |
|---|---|---|---|---|---|
| | Increase (MW) | Decrease (MW) | | | |
| G-5 | 45.50 | 0 | 0 | / | 15.60 |
| G-8 | 0 | 0 | 30.00 | / | 0 |
| G-11 | 32.07 | 0 | 0 | / | 0 |
| G-13 | 0 | 6.31 | 0 | 1.50 | 14.40 |
| C-30 | / | / | / | 1.50 | / |

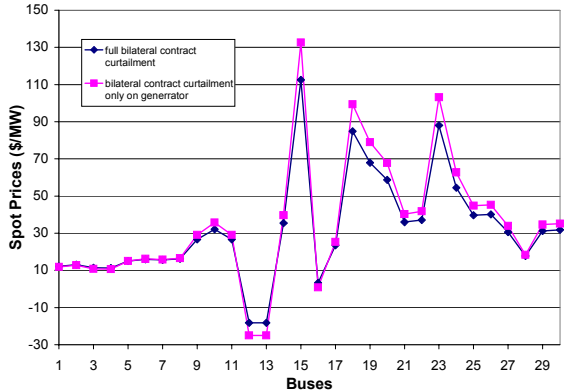

Fig.5. Spot prices for congestion management with bilateral contract curtailment

To demonstrate it, we change the curtailment strategy of bilateral contract to only decreasing the generation while keeping the load supplied by other available resources in RBM. The comparison of spot prices between these two situations is given in Fig.5. Compared with Fig.4 the fluctuation of spot prices in Fig.5 is very large, which means the congestion in this case is very serious. The price at bus 15 is prohibitively high because the consumer at bus 15 is the main cause of this congestion. Another phenomenon that should be noticed in Fig. 5 is that the spot prices at buses 12 and 13 are negative. Negative price at bus-13 means that decremental or bilateral curtailment bids from G-13 must have been accepted by the ISO to eliminate the network congestion. Negative price at bus-12 implies that L-12 can get payment in RBM from increasing its load since this action could be helpful to alleviate the congestion.

**Case 3**: Voltage limits violation. Change the active load at bus 26 from 3.5 MW to 33.5 MW. As a result, the voltage at bus 26 decreases to 0.82 p.u. and violates the lower limit 0.85 p.u. Reactive congestion management is run to eliminate the voltage violation. Fig. 6 shows the nodal prices of reactive power, which are the shadow prices of the nodal reactive power balancing equality constraints of the Q sub-problem. It is obvious that the reactive power price at bus 6 is prohibitively high. This penalty is reasonable, because the voltage violation is caused by the heavy active load at bus 26. This price signal also implies a planning requirement that some reactive power compensation devices should be installed at or near this spot.

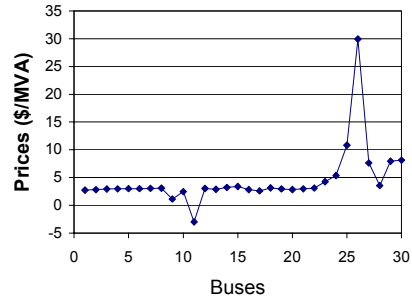

Fig.10 Effects of voltage violation on nodal prices

## VIII. CONCLUSIONS

In this article, a coordinated real-time optimal dispatch method for unbundled electricity markets has been proposed. The main features of this method are the following:

Three types of electricity contracts, which are bilateral energy contract, Pool energy auction contract and ancillary services contract, have been taken into account in the coordinated dispatch;

Balancing mechanism plays the key role in the proposed framework, where the ISO can meet the system imbalance and mitigate the network congestion by using various bids;

The adjustment of bus injection has been divided into several independent control variables according to unbundled contracts to embed all the possible bids in RBM into the objective of active power optimization;

The curtailment strategies of bilateral contracts have been integrated into the proposed method;

Test results demonstrate that the proposed coordinated dispatch method implemented with a modified OPF can deal with the system imbalance and network congestion simultaneously and successfully.

[3] WANG, X. and SONG, Y.H.: 'Advanced Real-Time Congestion Management through Both Pool Balancing Market and Bilateral Market', *IEEE Power Engineering Review*, Feb 2000, 20(2), pp.47-49

[4] WANG, X. and SONG, Y.H., and Lu, Q.: 'Apply Primal-dual Interior Point Linear Programming to Real-time Active and Reactive Power Congestion Management through Both Pool Balancing Market and Bilateral Market', 2000 IEEE/PES Winter Meeting, Paper No. WM2000-141, Singapore, Jan. 2000

[5] EI-KEIB, A.A., MA, X.: "Calculating Short-Run Marginal Costs of Active and Reactive Power Production", *IEEE Transactions on Power Systems*, Vol. 12, No. 2, May 1997, pp. 559-565
7